\documentclass[12pt, letterpaper, ]{article}

\usepackage[utf8]{inputenc}
\usepackage{amssymb}
\usepackage{multicol}
\usepackage[margin=1.98cm]{geometry}
\usepackage{setspace}
\usepackage{graphicx}
\usepackage{url}
\usepackage{color}
\usepackage[hidelinks]{hyperref}
\graphicspath{ {images/} }
\usepackage{tikz}
\usepackage{amsmath}

\begin{document}
	
\title{Peacock's Principle as a Conservative Strategy \bigskip \bigskip}

\author{\bigskip Iulian D. Toader\thanks{Institute Vienna Circle, University of Vienna, iulian.danut.toader@univie.ac.at}}

\date{\bigskip (forthcoming in \textit{Archive for History of Exact Sciences})}
	
\maketitle
	
\doublespacing

\begin{quote} 

\textbf{Abstract:} The view that Peacock's principle of permanence has been invalidated by Hamilton's introduction of non-commutative algebras has always seemed rather odd, in light of Peacock's favorable reception of quaternions and the endorsement of his principle by Hamilton. But the view is not just odd; it is incorrect. In order to show this, I critically analyze Peacock's attempts to reject possible exceptions to his principle, like the factorial function and an infinite series due to Euler. Then I argue that the principle of permanence is best understood as an expression of a conservative strategy, philosophically grounded in Hume's conception of the laws of reasoning, which advocates their preservation to the furthest extent possible, thus allowing exceptions, i.e., violations of these laws. On this reading, non-commutative multiplication does not invalidate Peacock's principle, if the reasons for violating commutativity outweigh the reasons for its preservation. Finally, I show that Hamilton followed a conservative strategy of precisely this sort when he developed his quaternionic calculus.

\end{quote}
	

\newpage 

\section{Introduction}
		
In a note included in his source book for the history of philosophy of mathematics, William Ewald wrote that ``Peacock's principle of the permanence of equivalent forms was a strait-jacket on his symbolical algebra, and led him to develop its rules in too close analogy with the laws of arithmetic.'' (Ewald 1996, 321) As such, symbolical algebra would not allow a geometrical interpretation of its operations, and more importantly, would leave no room for non-commutative operations. Indeed, in another note, Ewald wrote that Peacock's principle ``in effect guaranteed that symbolical algebra would always obey the familiar laws of arithmetic.'' (\textit{op. cit.}, 449) When characterizing Gregory's algebra, which was introduced about ten years after Peacock's, Ewald mentioned a significant step forward: ``Algebra has here been cut loose from the `principle of the permanence of equivalent forms' and from a too-close dependence on the arithmetic of the positive real numbers.'' (\textit{op. cit.}, 322) Despite allowing a geometrical interpretation of symbolical operations, this algebra still included no non-commutative operations --- cut loose from Peacock's principle, as it were, but not yet entirely. Not until a few years later, when Hamilton discovered the quaternions, and Cayley and Graves independently introduced non-associative operations on octonions.

The claim that Peacock's algebra would not allow a geometrical interpretation of its symbolical operations appears, however, to be based on a misunderstanding. In fact, his view was that an arithmetical interpretation of symbolical algebra is not necessary, but only necessarily possible, which obviously allows for many different interpretations, including a geometrical one (Toader 2024, 138). Peacock did regard the arithmetical interpretation as more convenient than a geometrical one (when convenience is determined by usefulness in applications) but later he came to see the geometrical interpretation (or more exactly, the goniometrical interpretation) as truly valuable. In any case, the focus of my paper is on the claim that, due to its commitment to the principle of the permanence of equivalent forms, Peacock's symbolical algebra would leave no room for non-commutative operations. This claim reiterates the old and common criticism that his principle has been invalidated by the development of increasingly abstract mathematical structures. The goal of my paper is to show that this criticism is mistaken.

The criticism is at least as old as Russell's fourth book, \textit{The Principles of Mathematics}, written in 1900 and published three years later. After mentioning the ``freer spirit towards ordinary Algebra'' that he thought guided the work of Hamilton, and ``the far
wider and more interesting problems of Universal Algebra'', as developed by Whitehead, Russell remarked: 

\begin{quote} \singlespacing

The possibility of a deductive Universal Algebra is often based upon a supposed principle of the Permanence of Form. Thus it is said, for example, that complex numbers must, in virtue of this principle, obey the same laws of addition and multiplication, as real numbers obey. But as a matter of fact there is no such principle. ... The principle of the Permanence of Form, therefore, must be regarded as simply a mistake: other operations than arithmetical addition may have some or all of its formal properties, but operations can easily be suggested which lack some or all of these properties. (Russell 1903, 377)

\end{quote} 

Of course, if one assumes that the principle of permanence stipulates the preservation of all arithmetical laws for all number systems, it is hardly surprising that one comes to believe, as Russell did, that non-commutative or non-associative operations are counterexamples to that principle. My attempt to refute this relentless criticism is motivated by the following brief observation. If the criticism were correct, it would make not only Peacock's favorable attitude towards quaternions (Crowe 1967, 34), but also Hamilton's application of the principle of permanence in his own mathematics (e.g., in the passage from quaternions to biquaternions; Fisch 2017, 201), if not disingenuous, then quite baffling. Furthermore, it would make the later adoption of the principle by 19th century mathematicians, especially Hankel, and by 20th century ones, like von Neumann (Toader 2025, ch. 3), rather mysterious. To ask the obvious, but seemingly neglected question, how could they have embraced the principle of permanence, if this had indeed been left behind, if mathematics had really been cut loose from it? Might violations of commutativity, as in quaternionic algebra, be compatible with Peacock's principle? If so, how?

To address these questions, the paper is structured in the following way. In section two, I present several formulations of Peacock's principle, which he offered between 1830 and 1845, and I note some of the salient differences between them. Then I explain his notion of equivalent forms, and his original conception of the relation between arithmetical algebra and symbolical algebra. I focus on his view that, in relation to symbolical algebra, a subordinate science like arithmetical algebra is suggestive, but not restrictive, in the sense that the meaning  of arithmetical operations is  not carried over to symbolical algebra. For comparison, I mention an alternative system of algebra, due to Kelland, which was emphatically advertised as a superior alternative to Peacock's system, for being both suggestive and restrictive. 

In section three, I begin addressing the problem of the universality of Peacock's principle. I revisit a long review by De Morgan, to discuss his objection that even if the equivalent forms in symbolical algebra satisfy the universality condition articulated by Peacock, i.e., even if they are to hold for any values of their variables, this would not guarantee the universality of his principle. Considering the differences between its formulations, I then show that the principle of permanence allows for the existence of equivalent forms that belong to symbolical algebra but not to arithmetical algebra. Thus, in at least one sense, this is not a universal principle.

In section four, I discuss the more serious objection that Peacock's principle is not universal because there are equivalent forms in arithmetical algebra that are not preserved in symbolical algebra. I focus on the factorial function, which clearly violates the universality condition. As Anna Bellomo recently noted, Peacock nevertheless attempted to show that the factorial function is preserved in symbolical algebra, since the factorial function can be derived from the Gamma function (Bellomo 2025, 158). But I argue that his attempt failed, because the Gamma function violates an invariance condition, which was articulated by Peacock himself in his wrestling with an infinite series that Euler had presented as a counterexample to the principle of permanence. The equivalent forms of symbolical algebra, Peacock argued, must be invariant, i.e., they must not change their constitution for any values of their variables. However, both Euler's infinite series and the Gamma function change their constitution for some values of their variables. More importantly, I think that Peacock's attempt illustrates the sort of deliberation that is needed to determine whether equivalent forms are  exceptional or not.

In section five, prompted by Mic Detlefsen's emphasis on a qualification of Peacock's principle (Detlefsen 2005, 283\textit{sq}), I turn to what I believe were the philosophical roots of the principle: Hume's conception of the laws of reasoning. As I understand him, Hume maintained that the laws of reasoning must be preserved to the furthest extent possible, due to their usefulness in applications, especially for belief formation. This allows for exceptions: although practically unsustainable, it would not be impossible to violate the laws. Peacock's principle should, I submit, be interpreted as an expression of this conservative strategy: symbolical algebra must preserve the equivalent forms of arithmetical algebra to the furthest extent possible, due to their usefulness in applications, especially for calculation. But this, likewise, allows for exceptions: equivalent forms in arithmetical algebra can be violated in symbolical algebra and, although practically inconvenient, the latter might be based on alternative, even entirely arbitrary, assumptions. If this interpretation is correct, then the existence of equivalent forms that belong to arithmetical algebra but not to symbolical algebra does not necessarily invalidate Peacock's principle. 

In section six, I argue that this is indeed the case with the commutativity of multiplication: although this arithmetical law is certainly violated in quaternionic algebra, this violation does not invalidate the principle of permanence. In fact, Hamilton followed Peacock's principle, attempting to preserve all arithmetical laws, but realized that there are reasons for abandoning commutativity that outweigh the reasons for preserving it. Hamilton's own account of the development of his quaternionic calculus describes precisely the sort of deliberation that I take to be part of the conservative strategy expressed by Peacock's principle. A similar conclusion may hold for more general mathematical structures.

\section{The Principle of Permanence}

George Peacock (9 April 1791 --- 8 November 1858) studied at Trinity College, Cambridge, where he became a lecturer and co-founded the Analytical Society. He was ordained as a priest and elected a Fellow of the Royal Society several years before he published his \textit{Treatise on Algebra} in 1830, with a second edition in two volumes, published in 1842 and 1845. Peacock was an active reader for the Cambridge Philosophical Society, and invited to report to the British Association for the Advancement of Science. In 1837, he became Lowndean Professor of Astronomy at Cambridge, but moved away in 1839, when appointed Dean of Ely Cathedral, holding on to both of these positions until the end of his life. But Peacock was first and foremost a mathematical innovator, one whose works have been aptly described as ``ingenious attempts to hold on to the old while being forced to grope creatively toward new options.'' (Fisch 2017, 14) The methodology that he articulated to guide such attempts is expressed by his principle to the permanence of equivalence forms, stated in the first edition of the \textit{Treatise} in the following way:

\begin{quote} \singlespacing
	
Whatever form is Algebraically equivalent to another, when expressed in general symbols, must be true, whatever those symbols denote.
		
Conversely, if we discover an equivalent form in Arithmetical Algebra or any other subordinate science, when the symbols are general in form though specific in their nature, the same must be an equivalent form, when the symbols are general in their nature as well as in their form. (Peacock 1830, 104; italics removed)

\end{quote}
	
In the report to the British Association for the Advancement of Science, this two-part statement of the principle is slightly modified:

\begin{quote} \singlespacing
 
Direct proposition: Whatever form is algebraically equivalent to another when expressed in general symbols, must continue to be equivalent whatever those symbols denote.
		
Converse proposition: Whatever equivalent form is discoverable in arithmetical algebra considered as the science of suggestion, when their symbols are general in their form, though specific in their value, will continue to be an equivalent form when the symbols are general in their nature as well as in their form. (Peacock 1833, 198-199; italics removed)

\end{quote}

Finally, in the two volumes of the second edition of his \textit{Treatise}, Peacock altered this formulation in a more radical manner:

\begin{quote} \singlespacing

[A]ll the results of arithmetical algebra which are deduced by the application of its rules, and which are general in form, though particular in value, are results likewise of symbolical algebra, where they are general in value as well as in form. (Peacock 1842, vi-vii)

\smallskip
	
Whatever algebraical forms are equivalent, when the symbols are general in form but specific in value, will be equivalent likewise when the symbols are general in value as well as in form. (Peacock 1845, 59)
		
\end{quote}

There are a few notable differences between these formulations of the principle, but I shall focus on only two of them here. First, some formulations are more general than others. Secondly, the earlier versions include two propositions, a ``direct'' one and its ``converse'', whereas the later ones include only versions of the latter. Before I discuss these differences, here is briefly what Peacock seems to have understood by the notion of equivalent forms in arithmetical algebra, where symbols --- variables like $x$, $m$, and $n$ below --- denote positive integers, and in symbolical algebra, where they are allowed to range over any quantities. Considering his own examples, one can see that what he called equivalent forms are the laws or rules for operations, such as $x(m+n) = xm + xn$ and $x y = y x$ and $(1+x)^m (1+x)^n = (1+x)^{m+n}$, but also the results (that is, the theorems) obtained through the application of rules, like the binomial theorem:

\begin{equation*}
(1+x)^n = 1+\frac{n}{1}x+\frac{n \cdot (n-1)}{1 \cdot 2}x^2+\frac{n \cdot (n-1) \cdot (n-2)}{1 \cdot 2 \cdot 3}x^3+ . . .
\end{equation*}

\bigskip

All such forms are mathematical equations. Peacock referred to them as ``equivalent forms'', in the plural, because the expressions on their left-hand and right-hand sides are numerically equal, for any values of their variables. But he also used the singular, when referring to the expression on either side as an ``equivalent form'' of the other. Thus, $(1+x)^m (1+x)^n$ is a form that is equivalent to $(1+x)^{m+n}$, just like the infinite series

\begin{equation*}
	1+\frac{n}{1}x+\frac{n \cdot (n-1)}{1 \cdot 2}x^2+\frac{n \cdot (n-1) \cdot (n-2)}{1 \cdot 2 \cdot 3}x^3+ . . .
\end{equation*}

\bigskip

is an equivalent form of $(1+x)^n$. More importantly, on Peacock's view, a subordinate science like arithmetical algebra is suggestive, in relation to symbolical algebra, but not restrictive. This means that, in the passage from arithmetical algebra to symbolical algebra, the latter must preserve the equivalent forms of the former, but not necessarily their meaning. In particular, the operations in symbolical algebra are not arithmetical: for example, ``+'' does not denote arithmetical addition, but symbolical addition, and the latter reduces to the former only when an expression like $(1+x)^n$ is interpreted over the positive integers. As such, according to Peacock, the rules and theorems of symbolical algebra do not have the same meaning as their arithmetical counterparts.

Peacock's view that arithmetical algebra is suggestive, but not restrictive, was questioned by his contemporaries. Philip Kelland, Professor of Mathematics in Edinburgh, argued that arithmetical algebra should be suggestive, but it should also be restrictive, in the sense that symbolical algebra should preserve not only the equivalent forms of arithmetical algebra, but also their meaning. For, he argued, if arithmetical meaning were not carried over, the results of symbolical algebra would be meaningless, and so their consistency with the results of arithmetical algebra could not be guaranteed (Kelland 1843, 120). Symbolical algebra, Kelland wrote,

\begin{quote} \singlespacing
	
	admits any arbitrary definitions whatever, provided it can be shewn that they do not lead to results inconsistent with arithmetic. The necessity that such be the case, however, renders this science a very different one from that which Mr Peacock expounds. In his, no prior knowledge of the nature of the symbols, no perception of the mutual agreement of the operations with each other or with another science is possible --- in ours such knowledge and perception is indispensable. (\textit{op. cit.}, 110sq.) 
	
\end{quote}

Kelland seems to have thought that one could show that the consistency requirement has been satisfied only if one had \textit{prior} knowledge of the nature of variables and operations in symbolic algebra. This knowledge would presumably be based not on their interpretation over the particular domain of positive integers, for that would be \textit{posterior} knowledge. The \textit{prior} knowledge that he demanded could be obtained by preserving not only equivalent forms, but their meaning as well, i.e., by carrying the meaning of operations over from arithmetical algebra to symbolical algebra. This indicates that Kelland did not allow for equivalent forms in symbolical algebra that are not equivalent forms in arithmetical algebra. Indeed, he proposed an alternative system, which he thought would be superior to Peacock's system, but which might be quite fittingly characterized as constrained by a straight-jacket:

\begin{quote} \singlespacing
	
	The system of algebra which I have endeavoured very briefly to develope takes arithmetic as its basis, and proceeds upwards by extending the definitions, and admitting to the forms it introduces, properties derived from the subordinate science. Arithmetic is \textit{suggestive}, but it is more, it is \textit{restrictive}. (\textit{op. cit.}, 120) 
	
\end{quote}

I turn now to the problem of the universality of Peacock's principle, beginning to address the following question: are all equivalent forms in arithmetical algebra and in symbolical algebra permanent, in the sense stipulated by the principle, or only some of them? First, I ask whether the principle allows the existence of equivalent forms in symbolical algebra that are not equivalent forms in arithmetical algebra. My answer is yes, with respect to rules, but no, with respect to theorems. Afterwards, I ask the more difficult question if the principle can accommodate equivalent forms in arithmetical algebra that are not equivalent forms in symbolical algebra. My answer is, again, yes; nevertheless, this does not invalidate Peacock's principle.

\section{The Direct Proposition}

In a long review (in two parts, each almost 20 pages) of Peacock's \textit{Treatise}, De Morgan had also reflected on the relation between arithmetical algebra and symbolical algebra (De Morgan 1835). He agreed with Peacock that a subordinate science must not be restrictive, but insisted that it should not be suggestive, either. De Morgan proposed that symbolical algebra should be based instead on entirely arbitrary assumptions. His argument is essentially as follows: if a subordinate science were suggestive, then the results of symbolical algebra would admit of many different interpretations. But since each interpretation would validate the same theorems, this would motivate a shift from meaning to pure form, which would imply that any subordinate science is, in principle, dispensable for the development of symbolical algebra: ``If we were to recommend any alteration, it should be to abandon, in a great measure, the \textit{science of suggestion}.'' (\textit{op. cit.}, 310) In a great measure, although not entirely abandoned, because a science of suggestion would have advantages in mathematical education, which a large part of De Morgan's review was concerned with.

With regard to Peacock's principle, as formulated in 1830, in the first edition of the \textit{Treatise}, De Morgan noted that while the converse proposition appears clear and unproblematic, the direct proposition should be scrutinized more carefully:

\begin{quote} \singlespacing
	
	For it seems to us, that the first part of the principle is that on which all the definitions are built. Their universality, and consequently the universality of the deductions legitimately drawn from them, is made the peculiar character of [symbolical] algebra. ... Mr. Peacock constructs those fundamental assumptions with no other intention than to justify the [universal] use of the principle. (\textit{op. cit.}, 308)
	
	
\end{quote}

De Morgan thought that Peacock was right to assume that the rules that define symbolical operations, and the symbolical results derived from them, hold universally, i.e., for any values of their variables. This agreed with De Morgan's own proposal that symbolical algebra be based on entirely arbitrary assumptions. But he also observed that the assumption that equivalent forms in symbolical algebra are universal, although certainly necessary, is not enough to establish the universality of Peacock's principle. According to De Morgan, what would be further needed is a justification of the direct proposition, and by consequence, a justification of the requirement that all equivalent forms in symbolical algebra be preserved in arithmetical algebra.

However, Peacock allowed for the existence of equivalent forms in symbolical algebra that are not equivalent forms in arithmetical algebra. In order to see this, we have to closely consider the differences between his formulations of the principle, quoted at the beginning of the previous section. One of these differences concerns the specification of a subordinate science. As Peacock formulated the 1833 converse proposition, this is less general than its 1830 version, since it does not mention any other subordinate science than arithmetical algebra. This might have created the impression that arithmetical algebra was actually thought to be, or maybe even ought to be, the only science of suggestion for the development of Peacock's symbolical algebra. It might have also induced the view that the latter's operations and results can admit of no other interpretation but the arithmetical one. However, as I have already noted in the introduction, this was not the case: indeed, the 1830 version clearly allows for non-arithmetical sciences of suggestion. Moreover, both versions of the direct proposition allow that general symbols can denote non-arithmetical quantities. Further, the 1842 formulation of the principle is less general than its 1845 version, in the same way, since it specifically mentions arithmetical algebra as a subordinate science; but in a different way as well, since it pertains only to results (that is, theorems), rather than all equivalent forms. Perhaps the most significant difference is Peacock's dropping the direct proposition altogether in these later formulations of his principle. What could explain this move? Did he realize that there are equivalent forms of symbolical algebra that do not remain equivalent forms when their variables range over the positive integers? Did he come to see that there might be rules of symbolical operations that fail in arithmetical algebra? 

It is fairly plausible to believe that such failure would be enough for Peacock to drop the direct proposition. This explanation has been proposed by Menachem Fisch, who quite naturally assumed that Peacock must have understood the direct proposition to apply to the rules of symbolical algebra, as well as to the theorems derived from them. But Fisch further made the following correct observation: ``Throughout the [first edition of the] \textit{Treatise} the Direct Proposition is applied only to low-level results --- never at the level of fundamental principle or rule of operation.'' (Fisch 2017, 193) However, it seems to me that the reason why Peacock never applied the direct proposition to rules, but exclusively to theorems, is that he actually understood it to be applicable only to theorems, and never to rules. This follows from his own justification of the direct proposition, conveniently given immediately after the 1833 formulation of the principle: 

\begin{quote} \singlespacing
	
The direct proposition must be true, since the laws of combination of symbols by which such equivalent forms [i.e., the forms that the direct proposition pertains to] are deduced, have no reference to the specific values of the symbols. (Peacock 1833, 199)

\end{quote}

This justification strongly indicates that the intended application of the 1833 direct proposition is to be restricted to theorems, i.e., equivalent forms that are \textit{deduced} in symbolical algebra via its rules of operation or laws of combination. Peacock thought that the direct proposition, restricted in this way to the theorems of symbolical algebra, was justified precisely because the rules of symbolical algebra are expressed in general symbols, without denotation of any particular objects, and they hold universally, for any values of their variables. Clearly, he wanted both universal rules and applicable theorems. But to ensure the applicability of the theorems of symbolical algebra, derived as they are from universal rules, he demanded that these theorems be transferred to arithmetical algebra (or any other subordinate science). This demand is, on my reading, expressed by the direct proposition.

But if applicable to theorems, and not to rules, then why did Peacock drop the direct proposition in the second edition of the \textit{Treatise}? After all, in the first volume of that edition, from 1842, his formulation of the principle explicitly pertains to theorems exclusively. It would have made complete sense to keep the direct proposition, similarly restricted, rather than dropping it. Did he realize that there might also be theorems of symbolical algebra that do not  hold for the positive integers? On closer analysis, the explanation of Peacock's move is that the direct proposition, applicable to theorems, but not to rules, was already assumed in his justification of an assumption behind the converse proposition:

\begin{quote} \singlespacing

[W]e may assume the existence of such an equivalent form in symbols which are general both in their form and in their nature, since it will satisfy the only condition to which all such forms are subject, which is, that of perfect coincidence with the results of arithmetical algebra, as far as such results are common to both sciences. (\textit{op. cit.}, 199)

\end{quote}

I think that this passage is to be understood as follows. To say that results (that is, theorems) are common to both sciences implies that they can be transferred from symbolical algebra, where their variables range over any objects, to arithmetical algebra, where they range over positive integers only. Such a transfer requires an arithmetical interpretation of the variables and the operations in symbolical algebra. Thus, we are justified to assume the existence of symbolical forms if this condition is satisfied: when transferred to arithmetical algebra, all the theorems of symbolical algebra must hold for positive integers. Peacock emphasized this condition repeatedly: ``symbolical conclusions [are] necessarily coincident with those of arithmetical algebra'' (\textit{op. cit.}, 201) But this expresses the very requirement stated by the direct proposition in its 1833 formulation, provided one acknowledges its restricted applicability, to theorems but not to rules. Peacock's condition of perfect coincidence might, however, be stronger than necessary. The weaker requirement that no theorem of symbolical algebra, when arithmetically interpreted, should contradict any arithmetical theorem, which was emphasized by Kelland (see the previous section), might have been enough to justify the existence of symbolical forms.

My explanation of Peacock's move, his dropping of the direct proposition in later formulations of the principle of permanence, goes some way towards challenging the claim that the principle was a straight-jacket, i.e., a universal principle stipulating that all equivalent forms in symbolical algebra must be equivalent forms in arithmetical algebra. It was nothing of this sort. For it turns out that, properly understood, the principle actually allows for rules of symbolical algebra that are not rules of arithmetical algebra, without allowing for theorems of symbolical algebra that are not (when interpreted over the positive integers) theorems of arithmetical algebra as well. But does it also allow for equivalent forms in arithmetical algebra --- either rules or theorems --- that are not equivalent forms in symbolical algebra? 

Since the affirmative answer to this question has been typically taken to invalidate the principle, I want to discuss Peacock's own attempts to reject two exceptions to the converse proposition: the factorial function and an infinite series due to Euler. Even if they were successful, these two attempts would not provide enough reasons for rejecting all possible exceptions. However, they help clarify additional conditions that Peacock thought must be satisfied by all equivalent forms that fall within the scope of his principle. We will see that when these conditions are taken into account, neither the factorial function, nor Euler's infinite series can be preserved in symbolical algebra.

\section{The Converse Proposition}

Consider first the factorial function, $n! = 1 \cdot 2 \cdot 3 \cdot ... \cdot (n - 2) \cdot (n - 1) \cdot n$. Although this holds for positive integers, it does not hold for all values of \textit{n}, since it obviously does not hold for the real numbers. In other words, $1 \cdot 2 \cdot 3 \cdot ... \cdot (n - 2) \cdot (n - 1) \cdot n$ is an equivalent form of $n!$ in arithmetical algebra, but not in symbolical algebra. Thus, since it fails the universality condition on symbolical forms, it would be enough to challenge the universality of Peacock's principle. The factorial function is but one example of what he, following Euler, called ``\textit{inexplicable functions} ... which are apparently
restricted by their form to integral and positive values of one or
more of the general symbols which they involve'' (\textit{op. cit.}, 210). 

Why apparently? Because, as Bellomo has emphasized (Bellomo 2025, 158), Peacock noted, following Legendre, that the factorial function can be derived from what he called a ``transcendent'' (Peacock 1833, 211) function $\Gamma(r+1) = \int_0^1 dx (log \frac{1}{x})^r$. This function, which he, again following Legendre, called the Gamma function, is defined for positive integers, but also for other values of \textit{r}, like the real numbers. The factorial function can be easily derived by noting that $\Gamma(r+1) = r \Gamma(r)$, about which Peacock said that it ``admits of all values of \textit{r}'' (\textit{loc. cit.}). Indeed, he insisted that $\Gamma(r+1) = r  \Gamma(r)$ holds ``whatever be the value of \textit{r}, whether positive or negative, whole or fractional. Legendre has apparently limited this equation to positive values of \textit{r}'', but Peacock considered this ``a restriction which is obviously unnecessary.'' (\textit{op. cit.}, 211\textit{sq}) Thus, like Peacock, one might be inclined to believe that $ \Gamma (r+1) = 1 \cdot 2 \cdot 3 \cdot ... \cdot (r - 2) \cdot (r - 1) \cdot r $ covertly or implicitly preserves the factorial function in symbolical algebra. 

But Legendre was, of course, right: when $r < 0$, the Gamma function diverges. Nevertheless, Peacock maintained that for negative values of \textit{r}, the Gamma function only changes its constitution: these values, he wrote, provide ``indications of a change in the constitution of the function'' (\textit{op. cit.}, 215). He implied that this change does not prevent the preservation of the factorial function in symbolical algebra. But it turns out that it violates one of the necessary conditions that Peacock, himself, imposed on the application of the principle of permanence --- the \textit{invariance} condition. To see what this is, we need to turn to his discussion of whether Euler's infinite series can be preserved in the passage from arithmetical to symbolical algebra.

Euler's infinite series seems to have been a genuine obsession for Peacock, for he kept coming back to it: he rejected it in three different ways, first in 1830, then in 1833, and one last time in 1845. His arguments are worth quoting in full. Before I do so, I should note, for historical accuracy, that Euler's paper, to which Peacock referred in all his responses, was \textit{not} published in the Petersburg Acts for 1774. This journal, with the full title \textit{Acta Academiae Scientiarum Imperialis Petropolitanae}, would be published only later, from 1777 until 1782, as a sequel to \textit{Novi Commentarii Academiae Scientiarum Imperialis Petropolitanae}, which ran from 1750 to 1776, and where Euler's paper, which was titled ``Demonstratio theorematis Neutoniani de evolutione potestatum binomii pro casibus, quibus exponentes non sunt numeri integri'', had in fact been published in 1775.

Here is Peacock's first response to Euler:

\begin{quote} \singlespacing
	
The opinion expressed in the text is expressly opposed to the great authority of Euler, who, in a species of preface to a demonstration of the Binomial Theorem which is given in the Petersburg Acts for 1774, has denied the universality of the principle of the permanency of equivalent forms ... 
The exception to the truth and universality of this principle which he quotes, is found in the very remarkable series

\begin{equation*}
\frac{1-a^m}{1-a}+ \frac{(1-a^m)(1-a^{m-1})}{1-a^2}+\frac{(1-a^m)(1-a^{m-1})(1-a^{m-2})}{1-a^3}+\&c.     
\end{equation*}

the law of the series being sufficiently manifest from that of its three first terms: if $m$ be a whole  [positive] number, the sum of this series is $m$, which he says is not the case for other values: but he appears to have confounded together, as is very commonly the case, the algebraical and arithmetical sum of the series, the first of which only is involved in the principle in question. (Peacock 1830, 517)
		
\end{quote}

The response here is essentially that Euler, like many others, including Cauchy (Peacock 1933, 192), had been confused about the relation between arithmetical and symbolical algebra. As we have seen in section two, Peacock saw the subordinate science as suggestive, but not restrictive, in the sense that the operations of arithmetical algebra and their symbolical counterparts do not have the same meaning. However, it's not clear how Euler's alleged confusion bears on the question whether this infinite series is an admissible exception to the principle of permanence. Did Peacock imply that even though, when \( m \) is not a positive integer, the arithmetical sum of the series does not equal \( m \), the symbolical sum would? In other words, considering the equation

	\begin{equation*}
	m = \frac{1-a^m}{1-a}+ \frac{(1-a^m)(1-a^{m-1})}{1-a^2}+\frac{(1-a^m)(1-a^{m-1})(1-a^{m-2})}{1-a^3}+\&c.     
\end{equation*}

did Peacock take it to be the case that, in symbolical algebra, the expression on the right-hand side is an equivalent form of the expression on the left-hand side, even though \textit{m} is not a positive integer? 

In any case, what he reported that Euler had said is, of course, correct: when \( m \) is not a positive integer, then if \( |a| < 1 \), the terms \( a^k \) tend to zero as \( k \) increases. In this case, the series might not converge to \( m \), for the exact sum depends on the specific value of \( m \) and \( a \). If \( |a| \geq 1 \), the terms \( a^k \) do not tend to zero, and the series may diverge. For \( a = 1 \), the series is undefined because the denominators \( 1 - a^k \) become zero. Peacock noted this last detail a few years later, in his second response:

\begin{quote} \singlespacing
	
Euler, in the Petersburgh Acts for 1774, has denied the universality of this principle, and has adduced as an example of its failure the very remarkable series

\begin{equation*}
	\frac{1-a^m}{1-a}+ \frac{(1-a^m)(1-a^{m-1})}{1-a^2}+\frac{(1-a^m)(1-a^{m-1})(1-a^{m-2})}{1-a^3}+\&c.     
\end{equation*}

\smallskip 

which is equal to $m$, when $m$ is a whole [positive] number, but which is apparently not equal to $m$, for other values of \textit{m}, unless at the same time $a = 1$: the occurrence, however, of \textit{zero} as a factor of the $(m+ 1)^{th}$ and following terms in the first case, and the reduction of every term to the form $\frac{0}{0}$ in the second, would form the proper indications of \textit{a change in the constitution of the equivalent function} corresponding to these values of $m$ and $a$. (Peacock 1833, 224; my emphasis)

\end{quote}
	
Peacock's point is that for certain values of $m$ and $a$, when the denominators \( 1 - a^k \) become zero and every term in the series reduces to $\frac{0}{0}$, the expression on the right-hand side of the equation

\begin{equation*}
	m = \frac{1-a^m}{1-a}+ \frac{(1-a^m)(1-a^{m-1})}{1-a^2}+\frac{(1-a^m)(1-a^{m-1})(1-a^{m-2})}{1-a^3}+\&c.     
\end{equation*}

changes its constitution and, thus, cannot be taken anymore as an equivalent form of the expression on the left-hand side. This subverts Euler's purported counterexample to the principle of permanence because, as Peacock implied in this second response, the equivalent forms that fall within the scope of the principle of permanence must be invariant. Invariance is considered to be an important condition on the application of the principle: it says that all equivalent forms in symbolical algebra, which fall within the scope of the principle, must be not only universal, i.e., they must hold for any values of their variables, but also invariant, i.e., they must not change their constitution for any values of their variables.

After ten years, Peacock returned once again to Euler's infinite series to provide yet another response ---  the third one:

\begin{quote} \singlespacing

Euler  ... in his celebrated proof of the series for $(1+x)^n$ ... denied the universal application of a principle equivalent to that of the permanence of equivalent forms, which alone could make it valid: he produced, as a striking exception to its truth, the very remarkable series [...]

\begin{equation*}
	\frac{1-a^m}{1-a}+ \frac{(1-a^m)(1-a^{m-1})}{1-a^2}+\frac{(1-a^m)(1-a^{m-1})(1-a^{m-2})}{1-a^3}+\&c.     
\end{equation*}

whose sum is $m$, when $m$ is a whole number, but not so for other values.	

A little consideration, however, will be sufficient to shew that the principle of the permanence of equivalent forms is not applicable to such a case: for if $m$ be a whole number, as in Arithmetical Algebra, the connection between $m$ and its equivalent series in the identical equation
	
	\begin{equation*}
		m = \frac{1-a^m}{1-a}+ \frac{(1-a^m)(1-a^{m-1})}{1-a^2}+\frac{(1-a^m)(1-a^{m-1})(1-a^{m-2})}{1-a^3}+\&c.     
	\end{equation*}
		
is not given, or, in other words, there is no statement or definition of the operation, by which we pass from $m$, on one side of the sign =, to a series under the specified form on the other, and there is consequently no basis for the extension of the conclusion to all values of the symbols, either by the principle of the permanence of equivalent forms or by any other: it is only when the results, which are general in form, but specific in value, are derived by processes which are definable and recognized, that they become the proper subjects for the application of this principle. (Peacock 1845, 452)

\end{quote}

The reason here adduced for rejecting Euler's infinite series, as falling outside the scope of the principle of permanence, is that it cannot be considered an equivalent form of \textit{m}, not even in arithmetical algebra. This is because no operation is definable or recognized that can lead from the expression on the left-hand side of the equation to that on its right-hand side. This might be taken to imply only that the infinite series does not exist necessarily, but may still be taken to exist as a hypothetical equivalent form of \textit{m}. Already in the preface to the first edition of his \textit{Treatise}, Peacock noted the following: ``When an equivalent form results from the performance of definable operations, its existence is \textit{necessary}."  (Peacock 1830, xviii). When it doesn't, he continued, it may still be said to exist hypothetically, ``only in virtue of the principle of the permanence of equivalence forms'' (\textit{loc. cit.}). Why wouldn't he then allow Euler's infinite series as a form that has a hypothetical existence in symbolical algebra? This is presumably because no operation, definable or not, which could lead from the expression on the left-hand side of the equation to that on its right-hand side, can be recognized.

To sum up, I take Peacock's view to have been that an equivalent form could be admitted as an exception to his principle only if, besides the universality of symbolical forms, it respects the following three conditions: the meaning of symbolical operations is not confounded with that of the arithmetical operations; equivalent forms must be invariant, i.e., their constitution should not depend on any particular values of their variables; and the operations on both sides of an equation should be definable and recognized. Clearly enough, Peacock thought that Euler's infinite series failed to satisfy each of these three conditions. For reasons already given above, the Gamma function fails at least the invariance condition. Thus, although in the latter case he may have wanted to conclude otherwise, the Gamma function cannot preserve the factorial function in the passage from arithmetical algebra to symbolical algebra. But then what should we say? Should we say that both Euler's infinite series and the Gamma function are inadmissible exceptions, that they fall outside the scope of the principle of permanence? Or should we say that even though they can be admitted as exceptions, they are not enough to invalidate Peacock's principle? If the latter, how might such exceptions be reconciled with it?

\section{Deliberative Conservatism}

One route towards reconciliation may be tracked to a qualification of the principle of permanence that was already acknowledged by Hankel (1867, 102\textit{sq}) and recently emphasized by Detlefsen: ``The Principle of Permanence says that the laws of basic arithmetic are to be retained \textit{to the fullest extent possible}.'' (Detlefsen 2005, 283sq; see also Toader 2021 for discussion.) Read as an instruction to preserve the equivalent forms of arithmetical algebra to the fullest or to the furthest extent possible in the development of symbolical algebra, Peacock's principle appears capable of accommodating exceptions. These would be equivalent forms, such as the factorial function and Euler's infinite series, that exist beyond the limits of the furthest possible extent to which equivalent forms are to be preserved. But how should one determine these limits? As Detlefsen duly put it: ``Hamilton's \textit{quaternions} ... and Graves's and Cayley's \textit{octonions}, ... both involve violations of basic arithmetic laws. But even the complexes involve violation of conditions that have some right to be regarded as standard arithmetic laws. ... What, then, should we say? Do the complexes violate the Principle of Permanence or do they not?'' (\textit{loc. cit.}) Do the quaternions, for that matter, and the octonions violate the principle or do they not? Where should the limits be set for the furthest possible extent to which equivalent forms are to be preserved? What might determine these limits?

All these questions can be answered by interpreting Peacock's principle in the light of Hume's conception about the permanence of the laws of reasoning. This interpretation makes sense of the above qualification of the principle, for Hume, as I read him, maintained that the laws of reasoning must be preserved to the furthest extent possible, and for mostly practical reasons rather than, say, as transcendental conditions for knowledge. This allows for exceptions, if they are backed by sufficiently strong reasons: although violation is practically unsustainable in general, it would not be impossible to violate the laws of reasoning. Interpreting the principle of permanence as an expression of a conservative strategy of this sort implies that the equivalent forms of arithmetical algebra must be preserved to the furthest extent possible. They must be preserved, as Peacock emphasized, because they are useful in applications, especially for calculations, rather than indispensable for the development of symbolical algebra. His principle can, likewise, allow for exceptions: equivalent forms could be abandoned when this is supported by sufficiently strong reasons. As I have noted already, Peacock admitted that symbolical algebra could be based on entirely arbitrary assumptions, although he thought that this would be practically inconvenient. But furthermore, this interpretation of the principle also clarifies how to determine the limits of the furthest extent possible to which equivalent forms are to be preserved. The determination involves deliberation, i.e., weighing the reasons for abandoning equivalent forms against the reasons for preserving them. 

My point here is not primarily a historical one: I do not claim that the way Peacock understood his principle was actually influenced by his reading of Hume. Yet such a claim would hardly be far-fetched. For it is perfectly reasonable to assume that Peacock had read Hume. As a cleric and founder of the Analytical Society, he could not have ignored Hume's philosophical views. We know that Herschel was familiar with these views, since he engaged with them, especially on causation, in his 1830 \textit{Preliminary Discourse}. His father, William Herschel, a German immigrant and royal astronomer, who discovered Uranus, had known Hume personally. Babbage also discussed Hume's views, on miracles, in his 1837 \textit{Ninth Brigdewater Treatise}. When Peacock's first edition of the \textit{Treatise} appeared in 1830, Hume's \textit{Philosophical Works}, in four volumes, had just been published in 1826. Besides, Peacock was known as an adherent of the Whig Party (Pycior 1981, 26), which might not be insignificant in this context, given the widespread view that Hume was the real ancestor of the Whigs (Snyder 2006, 14). It is true that, given the predominant practice of citation at the time, Hume is never mentioned in any of Peacock's mathematical writings. Nevertheless, I think that Hume's conception of the permanence of the laws of reasoning constitutes the philosophical background against which Peacock's principle is best understood. Of course, this is not to say that different backgrounds might not be relevant and useful for understanding the principle, for instance, the early
Cambridge Philosophical Society culture of natural history (cf. Lambert 2013) or a Lockean evidentialist epistemology of mathematics (cf. Richards 1992), but a detailed analysis of such proposals must be deferred elsewhere.

To recall the relevant details of Hume's conception of the laws of reasoning, let us consider a famously controversial passage from his \textit{Treatise} (in section 1.4.4.1), where he provided a justification for the rejection of what he called ``ancient philosophy'':

\begin{quote} \singlespacing

In order to justify myself, I must distinguish in the imagination betwixt the principles which are permanent, irresistible, and universal; such as the customary transition from causes to effects, and from effects to causes: And the principles, which are changeable, weak, and irregular; [...] The former are the foundation of all our thoughts and actions, so that upon their removal human nature must immediately perish and go to ruin. The latter are neither unavoidable to mankind, nor necessary, or so much as useful in the conduct of life; but on the contrary are observ’d only to take place in weak minds, and being opposite to the other principles of custom and reasoning, may easily be subverted by a due contrast and opposition. For this reason the former are received by philosophy, and the latter rejected. (Hume 1739, 148)

\end{quote}
		
The principles that Hume said are rejected by modern philosophy are \textit{ad hoc} rules of reasoning based, for instance, on an inclination to attribute occult qualities to material objects or on the fear to be tormented by spectres. Such rules are considered changeable since they can be easily abandoned in the face of evidence and normal evidence-assessing abilities, and they lack practical utility anyway, even if occult qualities and spectres are understood as fictions. The principles that Hume said are adopted by modern philosophy are the laws of probable (or inductive) and demonstrative (or deductive) reasoning. They are considered permanent, but this does not mean that it is impossible to abandon them. He explicitly envisaged their removal, and although he certainly implied that they should not be abandoned, he did think that they could be abandoned, as indeed they sometimes actually are.

However, even though the laws of reasoning could be abandoned, they could not be abandoned easily, or not as easily as the \textit{ad hoc} rules. This is because abandoning the laws would be practically inconvenient, and sometimes even unsustainable. Indeed, it may cause otherwise avoidable death. Imagine someone on a tall building considering the conflict between the customary transition from causes to effects and the idea that the ground below has the quality of kindly embracing any falling body. In most cases, the former subverts the latter, or perhaps more accurately, since principles are conceived of as belief-forming mechanisms (Loeb 2010, 83), the belief that one will fall to one's death easily subverts the belief that one will be kindly embraced upon falling. But it is not impossible for subversion to go, tragically, in the opposite direction. The point is that, in such cases, there are reasons for abandoning the laws of inductive reasoning that outweigh the reasons for preserving them. 

The same point may be perhaps better illustrated in relation to the laws of deductive reasoning. Consider, for instance, \textit{ex falso quodlibet}. Although, on Hume's conception, one should preserve this law, there are situations in which one could have indeed strong reasons to reject it. One such rejection, which Hume might have well been aware of, had been hinted at by Jonathan Swift in one of his fabulous satiric essays:

\begin{quote} \singlespacing
	
For, of what use is freedom of thought if it will not produce freedom of action, which is the sole end, how remote soever in appearance, of all objections against Christianity? And therefore the free-thinkers consider it as a sort of edifice, wherein all the parts have such a mutual dependence on each other that if you happen to pull out one single nail, the whole fabric must fall to the ground. This was happily expressed by him who had heard of a text brought for proof of the Trinity, which in an ancient manuscript was differently read; he thereupon immediately took the hint, and by a sudden deduction of a long sorites most logically concluded ``Why, if it be as you say, I may safely whore and drink on, and defy the parson.'' (Swift 1708, 226)

\end{quote}

The ironic characterization of the free-thinkers' logical inference as ``a sudden deduction of a long sorites'' is remarkable. But the point is that Swift implied here that the reasons for preserving any logical law governing this inference, such as \textit{ex falso quodlibet}, are outweighed by the reasons for abandoning it. For his belief that the behavior he described is immoral subverts the belief that such behavior can be justified by perceived inconsistencies in the doctrine of the Trinity. More generally, Hume allowed that the laws of reasoning could be subverted and abandoned, as indeed they sometimes are, in exceptional cases. This underlies a conservative strategy, which stipulates not only that the laws should be preserved to the furthest extent possible, but also that preservation depends on the weighing of relevant reasons. And when interpreted as an expression of this conservative strategy, as I think it should be interpreted, Peacock's principle likewise allows for the possibility that the equivalent forms of arithmetical algebra be subverted and abandoned when this is the outcome of deliberation, of a considered weighing of relevant reasons. 

Now, if one accepts this interpretation of Peacock's principle, then it is not difficult to see that the principle can accommodate exceptions such as the factorial function and Euler's infinite series, discussed in the previous section of this paper. But as I show in the next section, it is also not difficult to see that it can accommodate as well Hamilton's rejection of the law of commutativity in quaternionic algebra. Rejection is compatible with the stipulation that equivalent forms should be preserved to the furthest extent possible, if there are reasons for abandoning equivalent forms that outweigh the reasons for keeping them. Before I offer an explicit argument for this compatibility claim, illustrated by the case of quaternions, let me first point out that the widespread view, according to which Hamilton's development of his algebra can be considered as having invalidated Peacock's principle, actually neglects the available historical evidence.

\section{The Quaternions}

In a book-length presentation of quaternions (Hamilton 1853), comprising his lectures given as a Professor of Astronomy at Trinity College in Dublin, Hamilton, knighted in 1835 for his work in optics, described what he called a new mathematical method, i.e., ``the Method or Calculus of Quaternions'', the basis of which was his view of mathematics as a science of temporal and spatial order. This view led him to introduce a version of vectorial calculus, which he called symbolical geometry, and in which he ``sought to imitate the Symbolical Algebra of Dr. Peacock'' (\textit{op. cit.}, 61). Indeed, Hamilton considered the ordinal relations between moments and steps in time, and those between points and lines in space, as suggestive for the operations with ``directed distances'' (i.e., vectors) in symbolical geometry, but not restrictive, in exactly Peacock's sense discussed above, in section two. Furthermore, in passing from symbolical geometry to quaternionic algebra, the rules for operations with vectors are taken to suggest the rules for purely algebraical operations with quaternions, in the same way in which Peacock had taken the rules of arithmetical algebra to suggest those of symbolical algebra. The indebtedness to Peacock is unmistakable and properly acknowledged:  

\begin{quote} \singlespacing

Thus, in the phraseology of Dr. Peacock, we should have a very wide ``science of suggestion'' (or rather, suggestive science) as our \textit{basis}, on which to build up afterwards a new structure of purely \textit{symbolical generalization}, if the \textit{science of time} were adopted, instead of merely Arithmetic, or (primarily) the doctrine of \textit{integer number}. [...] And in the passage which I have made (in the Seventh Lecture), from \textit{quaternions} considered as \textit{real} (or as geometrically \textit{interpreted}), to \textit{biquaternions} considered as \textit{imaginary} (or as geometrically \textit{uninterpreted}), but as symbolically \textit{suggested} by the generalization of quaternion formulae, it will be perceived, by those who shall do me the honour to read this work with attention, that I have employed a \textit{method of transition}, from \textit{theorems proved} for the \textit{particular} to \textit{expressions assumed} for the	\textit{general}, which bears a very close \textit{analogy} to the methods of Ohm and Peacock. (\textit{op. cit.}, 16)
\end{quote}

However, my claim that the quaternions cannot be taken to invalidate Peacock's principle is supported not only by textual evidence that Hamilton applied it in his mathematics. After all, and not unusually at all, maybe he just paid lip service to a senior colleague. Still, when one understands what actually led him to reject commutativity, one realizes that Hamilton carefully applied the conservative strategy described above: he thought that this law should be preserved, but dropped it after weighing the reasons for preserving it against the reasons for abandoning it. 

In the preface to his book, Hamilton recalled a series of early attempts to extend the operation of multiplication of lines from two to three spatial dimensions. The first geometrical construction that he considered, in 1831, did not satisfy him because ``it did not preserve the \textit{distributive principle} of multiplication; a \textit{product of sums} not being equal, in it, to the \textit{sum of the products}.'' (\textit{op. cit.}, 36) He then mentioned ``another construction, of a somewhat similar character, and liable to similar objections, for the product of two lines in space'', which he said occurred to him in 1835 and a year later, independently, to his friend John T. Graves. But Hamilton thought this defective, too, for the same reason, as he found it ``to violate the distributive principle''. A third construction, proposed by John's younger brother Charles, assumed a representation of lines by triplets, and constructed the multiplication of two lines as the multiplication of two triplets (or ternions), $x+iy+jz$ and $x'+iy'+jz'$, where $i$ and $j$ are imaginary units, i.e., $i^2=j^2=-1$, such that a third triplet $x''+iy''+jz''$ would be the product. This operation was not distributive, either, but it was commutative. 

Inconsistency with the law of distributivity was, thus, the stated reason for abandoning these early constructions, which indicates that Hamilton thought that the usual laws of multiplication, including not only distributivity, but associativity and commutativity as well, should be preserved, assuming of course that they could be preserved. Indeed, when he came to reconsider the question, he said so explicitly:

\begin{quote} \singlespacing
	With such preparations as I have described, I resumed
	(in	1843) the endeavour to adapt the general conception of triplets to the multiplication of lines in space, resolving to \textit{retain} the \textit{distributive} principle, with which some formerly conjectured systems had been inconsistent, and at first supposing that I \textit{could}
	preserve the \textit{commutative} principle also. (\textit{op. cit.}, (43))
	
\end{quote}

Supposing at first that commutativity holds for the multiplication of any two lines, and thus that it holds for the imaginary units as well, so that $ij = ji$, Hamilton realized that there is no satisfactory way to represent the line product as a triplet. This gave him one reason to doubt that the law of commutativity could be preserved. As he further explained in more detail in the lectures, he also found that the multiplication of lines would not yield unique results, which gave him an even stronger reason to doubt commutativity: 

\begin{quote} \singlespacing
	
The multiplication of \textit{lines} among themselves has been shewn to give \textit{different results}, according as the factors have been taken in one or in another \textit{order}; from which it follows, by still stronger reason, that the \textit{multiplication of quaternions} is \textit{not} generally a \textit{commutative} operation. (\textit{op. cit.}, 132)
	
\end{quote}

Nevertheless, instead of abandoning the law of commutativity right away, Hamilton still tried to preserve it. In order to do so, he replaced the line representation proposed by Charles Graves with an alternative representation, as the triplet $ix+jy+kz$, and was led to a different definition of multiplication, which had $ij = +k$ and $ji = -k$, with $i^2 = j^2 = -1$. However, this left $k$ undetermined. Subsequently, Hamilton settled on the now well-known representation of a line as a quaternion $x+iy+jz+kw$, which takes the multiplication of quaternionic units to be non-commutative: $i^2 = j^2 = k^2 = -1$, $ij = -ji = k$, $jk = -kj = i$, and $ ki= - ik = j$, and then proved that commutativity fails for the product of any two ``rectangular vectors'' (i.e., mutually perpendicular lines): $\alpha\beta = - \beta\alpha$. He frankly acknowledged that this ``might at first seem strange'' (\textit{op. cit.}, 48) but noted that ``we are compelled, by considerations which appear more primary, to \textit{give up the commutative property} of multiplication, as not holding \textit{generally} for \textit{lines}.'' (\textit{op. cit.}, 51) Further, Hamilton showed that, unlike commutativity, the laws of distributivity and associativity are preserved in quaternionic algebra, and moreover, he argued that they should both be preserved due to their usefulness in applications:

\begin{quote} \singlespacing

The Calculus of Quaternions would ... be extremely incomplete, if it were \textit{permanently} deprived of the use of either of these two important principles: and indeed the \textit{combination of both} is \textit{essential} in many of its more advanced applications. (\textit{op. cit.}, 494)
	
\end{quote}

I think that this description of Hamilton's own thinking process, which led him eventually to the introduction of quaternions, fully supports my claim that the violation of the arithmetical law of commutativity in quaternionic algebra cannot be taken to invalidate Peacock's principle. As we have seen, not only did Hamilton apply the principle of permanence in developing his calculus, but he saw no incompatibility at all between this application and his rejection of commutativity. This is not surprising, if we understand Peacock's principle as an expression of the deliberative conservatism that I have described in the previous section. This interpretation may also explain Peacock's approval of Hamilton's work and dispel the mystery of the principle's later adoption by other mathematicians.

\section{Conclusion}

In this paper, I reconsidered the criticism according to which Peacock's principle of the permanence of equivalent forms has been invalidated by the development of mathematics, and in particular by Hamilton's introduction of quaternions. This criticism, raised by Russell, has been repeatedly made against this principle, fairly recently by Ewald (see also Gray 2008, 28 and Bellomo 2025, 167). I argued that the criticism is mistaken. 

After clarifying Peacock's notion of equivalent forms and his conception of the relation between arithmetical algebra and symbolical algebra, I discussed several differences between the various formulations of the principle, which he gave between 1830 and 1845. Then I argued that these differences indicate that Peacock allowed the existence of symbolical rules that do not belong to arithmetical algebra. Further, I critically discussed his attempts to preserve the factorial function and to reject an infinite series due to Euler, but I argued that these are exceptions that can be accommodated, if Peacock's principle is interpreted as an expression of a conservative strategy of the sort that Hume articulated with respect to the laws of reasoning. On this interpretation, the principle  stipulates the preservation of equivalent forms to the furthest extent possible, and thus can allow that equivalent forms can be subverted and abandoned, if there are reasons for doing so that outweigh the reasons that support their preservation. 

Finally, I showed that this is precisely the case with Hamilton's quaternions. Following the steps that led him to abandon commutativity, in the passage from symbolical geometry to the calculus of quaternions and biquaternions, it is clear that he thought that all the usual laws, including the laws of distributivity, associativity, and commutativity, should be preserved. But he also realized that the reasons for abandoning the commutativity of multiplication outweigh the reasons for preserving it. In contrast, he maintained that the reasons for dropping associativity and distributivity are outweighed by the reasons for keeping them. Other mathematicians came to differ about this, but I expect that my interpretation of Peacock's principle can withstand potential criticisms analogous to the one discussed above. Yet others have reinterpreted the principle and defended its universal validity at the cost of rejecting quaternions as numbers (Schubert 1894, 578). But an account of this reinterpretation belongs elsewhere.









\section{Acknowledgments}

Thanks to Jeremy Gray, Iustin P. Toader, and an anonymous referee for comments that helped improve the paper. The Austrian Science Fund (FWF) is gratefully acknowledged for financial support through 
Grant 10.55776/PAT3440123.

\section{References} \onehalfspacing
		
Bellomo, A. (2025) Peacock's Principle of Permanence and Hankel's Reception, \textit{HOPOS:} 

\textit{The Journal of the International Society for the History of Philosophy of Science}, 15, 150--176

\smallskip

\noindent Crowe, M. J. (1967) \textit{A History of Vector Analysis. The Evolution of the Idea of a Vectorial System}, 

The University of Notre Dame Press, 2nd ed., 1985, Dover Publications 

\smallskip

\noindent De Morgan, A. (1835) \textit{A Treatise on Algebra}, by George Peacock, \textit{The Quarterly Journal of} 

\textit{Education}, 9, 293--311
		
\smallskip

\noindent Detlefsen, M. (2005) Formalism. \textit{The Oxford Handbook of Mathematics and Logic}, ed. by 

S. Shapiro, Oxford University Press, 236--317
				
\smallskip

\noindent Ewald, W. B. (1996) \textit{From Kant to Hilbert: A Source Book in the Foundations of Mathematics}, 

vol. I, Oxford University Press 
		
\smallskip

\noindent Fisch, M. (2017) \textit{Creatively Undecided. Toward a History and Philosophy of Scientific Agency}, 

The University of Chicago Press

\smallskip

\noindent Gray, J. (2008) \textit{Plato's Ghost: The Modernist Transformation of Mathematics}, Princeton 

University Press

\smallskip

\noindent Hamilton., W. (1853) \textit{Lectures on Quaternions}, Dublin: Hodges and Smith

\smallskip

\noindent Hankel, H. (1867) \textit{Vorlesungen über die komplexen Zahlen und ihre Funktionen}. Leipzig: Voss

\smallskip

\noindent Hume, D. (1739) \textit{A Treatise of Human Nature}, ed. 2007, Clarendon Press
		
\smallskip

\noindent Kelland, Ph. (1843) \textit{Lectures on the Principles of Demonstrative Mathematics}. Edinburgh:

Adam and Charles Black

\smallskip

\noindent Lambert, K. (2013) A Natural History of Mathematics: George Peacock and the Making of English 

Algebra, \textit{Isis}, 104, 278–302

\smallskip

\noindent Loeb, E. L. (2010) \textit{Reflection and the Stability of Belief. Essays on Descartes, Hume, and Reid}, 

Oxford University Press

\smallskip

\noindent Peacock, G. (1830) \textit{A Treatise on Algebra}, Cambridge: J. \& J. J. Deighton

\smallskip

\noindent Peacock, G. (1833) Report on the Recent Progress and Present State of Certain Branches of

Analysis. \textit{Report to the Third Meeting of the British Association for the Advancement of} 

\textit{Science}, London, 3, 185--352 

\smallskip

\noindent Peacock, G. (1842) \textit{A Treatise on Algebra}. 2nd ed., vol. I, Cambridge: J. \& J. J. Deighton

\smallskip

\noindent Peacock, G. (1845) \textit{A Treatise on Algebra}. 2nd ed., vol. II, Cambridge: J. \& J. J. Deighton

\smallskip

\noindent Pycior, H. M. (1981) George Peacock and the British Origins of Symbolical Algebra. \textit{Historia} 

\textit{Mathematica}, 8, 23--45

\smallskip

\noindent Richards, J. (1992) God, Truth, and Mathematics in Nineteenth-Century England. In M. J. Nye, 

J. L. Richards, and R. H. Stuewer (eds.) \textit{The Invention of Physical Science. Intersections of} 

\textit{Mathematics, Theology and Natural Philosophy Since the Seventeenth Century}, Springer, 

Dordrecht, 51--78

\smallskip

\noindent Russell, B. (1903) \textit{The Principles of Mathematics}, Cambridge University Press

\smallskip

\noindent Schubert, H. (1894) Monism in Arithmetic, \textit{The Monist}, 4, 561--579

\smallskip

\noindent Snyder, L. J. (2006) \textit{Reforming Philosophy. A Victorian Debate on Science and Society}, 

The University of Chicago Press

\smallskip

\noindent Swift, J. (1708) An Argument to prove that the Abolishing of Christianity in
England, may as 

things now stand, be attended with some Inconveniences, and perhaps not produce those many 

good Effects proposed thereby. In \textit{Major Works}, ed. by A. Ross and D. Woolley, 1984, 

Oxford University Press, 217--227

\smallskip

\noindent Toader, I. D. (2021) Permanence as a principle of practice. \textit{Historia Mathematica}, 54, 77--94

\smallskip

\noindent Toader, I. D. (2024) Is Bohr’s correspondence principle just Hankel’s principle of permanence? 

\textit{Studies in History and Philosophy of Science}, 103, 137--145

\smallskip

\noindent Toader, I. D. (2025) \textit{Rules and Meaning in Quantum Mechanics}, European Studies in Philosophy 

of Science, Springer

\end{document}